\documentclass[a4paper]{amsart}

\makeatletter
\def\subsection{\@startsection{subsection}{2}%
  \z@{.9\linespacing\@plus.7\linespacing}{.3\linespacing}%
  {\normalfont\bfseries}} 

\def\subsubsection{\@startsection{subsubsection}{3}%
  \z@{.9\linespacing\@plus.7\linespacing}{.1\linespacing}%
  {\normalfont\itshape}} 

\def\paragraph{\@startsection{paragraph}{4}%
  \z@\z@{-\fontdimen2\font}%
  {\normalfont\bfseries}}
\makeatother

\usepackage[utf8]{inputenc}
\usepackage[ISO]{diffcoeff}
\usepackage{mathpazo,
	xfrac,
	amssymb,
	float,
	graphicx,
	xcolor,
	tabularx,
	booktabs,
	cite,
	changepage,
	multirow,
	fix-cm,
	microtype
	}
\usepackage[detect-all]{siunitx}
\usepackage[overload]{empheq}
\usepackage[foot]{amsaddr}

\newcommand{\eg}{\textit{e.\,g.}}
\newcommand*{\tran}{^{\mkern-1.5mu\mathsf{T}}}
\DeclareMathOperator{\diag}{diag}
\newcommand{\cg}{\textnormal{\textsl{g}}}

\usepackage[
	hmarginratio={1:1},% equal left and right margins
  	vmarginratio={1:1},% equal top and bottom margins
  	textwidth=410pt,% new text width
	textheight=680pt,
  	heightrounded,% always useful
	]{geometry}
\newlength{\extralength}
\newlength{\fulllength}
\definecolor{bluecite}{HTML}{0875b7}
\usepackage[colorlinks=true, 
	linkcolor=bluecite, 
	citecolor=bluecite, 
	urlcolor=bluecite]{hyperref}
\usepackage{doi}

\begin{document}

\title[]{Simpson variational integrator for nonlinear systems: a tutorial on the Lagrange top}

\author{Juan Antonio Rojas-Quintero\,$^{1,\star}$}
\address{$^1$\quad\normalfont Secihti/Tecnológico Nacional de México/I. T. Ensenada, Ensenada, 22780, B.C., Mexico}
\email{$^\star$jarojas@secihti.mx}

\author{François Dubois\,$^{2,3}$}
\address{$^2$\quad\normalfont Université Paris-Saclay, Laboratoire de Mathématiques d’Orsay, Orsay, 91400, France}
\address{$^3$\quad\normalfont Conservatoire National des Arts et Métiers, LMSSC, Paris, 75141, France}

\author{Frédéric Jourdan\,$^{4}$}
\address{$^4$\quad\normalfont IMT Atlantique, LS2N, UMR CNRS 6004, 44307 Nantes, France}

\date{\today, Ensenada, B.C., Mexico}

\begin{abstract}
This contribution presents an integration method based on the Simpson quadrature. The integrator is designed for finite-dimensional nonlinear mechanical systems that derive from variational principles. The action is discretized using quadratic finite elements interpolation of the state and Simpson's quadrature, leading to discrete motion equations. The scheme is implicit, symplectic, and fourth-order accurate. The proposed integrator is compared with the implicit midpoint variational integrator on two examples of systems with inseparable Hamiltonians. First, the example of the nonlinear double pendulum illustrates how the method can be applied to multibody systems. The analytical solution of the Lagrange top is then used as a reference to analyze accuracy, convergence, and precision of the numerical method. A reduced Lagrange top system is also proposed and solved with a classical fourth-order method. Its solution is compared with the Simpson solution of the complete system, and the convergence order of the difference between both is consistent with the order of the classical method.
\end{abstract}

\maketitle


\section{Introduction}

The evolution of systems in applied sciences is commonly characterized by differential equations for which exact solutions are rarely available. Numerical integration is thus required in most cases. Variational integrators are a class of numerical integrators that possess properties such as symplecticity and momentum conservation that lead to no artificial energy dissipation. Therefore, variational integrators have been widely used in certain disciplines. One of the first symplectic integrators was proposed by de Vogelaere \cite{vogelaere1956}, and later, the good conservation of the integrals of motion was deemed to be advantageous for long-term and large-scale simulations in celestial mechanics and molecular dynamics (see, \eg, \cite{gladman1991,gray1994}).

For Hamiltonian and Lagrangian systems, both the framework of variational error analysis and the exact discrete Lagrangian have been established by Marsden and West in \cite{marsden2001}, which also introduced variational integrators of the Galerkin type (later summarized in \cite{leok2005}). Variational integrators constructed from a finite-dimensional approximation space and a quadrature rule are quasi-optimal, meaning that their accuracy is related to the best approximation error of the approximation spaces being used \cite{hall2015}. On the practical side, it has been demonstrated that variational integrators can be constructed by following different paths. For example, Hairer et al. \cite{hairer2006} have shown that the very popular Runge-Kutta methods can be formulated such that symplecticity is ensured by using either collocation methods or partitioning schemes. Later, Leok et al. showed that variational integrators can be designed using any underlying one-step method \cite{leok2012} or by using Taylor expansions \cite{schmitt2017}. 

The systematic construction of variational integrators is a well-established technique. Consequently, there has been little research focusing on the construction of low-order integrators that could be of interest for the simulation of nonlinear control systems, for example. Linear systems are usually treated in the literature to illustrate the key properties of variational integrators. Nonlinear systems with separable Hamiltonians are also routinely simulated (see \eg\ \cite{mcLachlan1992,ober2014,blanes2017,amodio2022}). Splitting methods are particularly well-adapted to solve separable systems \cite{hairer2006,sanzSerna2018}. However, inseparable Hamiltonian systems are more complicated to integrate. One strategy involves transforming the inseparable system into a separable one using the extended phase space approach \cite{pihajoki2014}, or by combining that approach with symmetric projections \cite{hairer2006,jayawardana2022}, in the framework of Runge-Kutta methods. 

In general, the solution process of nonlinear systems involves Newton's method. However, implementing the available variational integrators on inseparable nonlinear systems with multiple degrees of freedom is not trivial. This contribution presents a fourth-order variational integrator formulated specifically for nonlinear systems. The proposed method is compatible with multibody systems characterized by inseparable Hamiltonians, deriving from Hamilton's principle \cite{rojo2018}. The three degrees of freedom Lagrange top system is taken to exemplify the usage of the proposed integrator. A \emph{Lagrange top} is any rotational axisymmetric rigid body with a fixed contact-point with a flat surface \cite{audin1999,cushman2015}, also known as \emph{heavy symmetrical top} \cite{goldstein2002}. Lagrange analyzed the dynamics of this body in \textit{Mécanique Analytique} \cite{lagrange2009}, hence the name of the object. The Lagrange top is usually solved using quaternions \cite{bobenko1999,wolf2014,lin2024}, which require a specific formulation of discrete Lagrangian mechanics on Lie groups \cite{bobenko1999,marsden1999}. However, Euler angles are used in this work to analyze two types of motion of the symmetric top: loops and cusps. Lie group formulations remain out of the scope of this particular work.

The proposed integrator is built by applying Simpson's quadrature to discretize the least action principle using the discrete Lagrangian as in \cite{dubois2023_GSI,dubois2023_FVCA}. Lagrange finite elements in the $\mathbb{P}_2$-space are used for interpolation. The resulting discrete Euler-Lagrange equations form an implicit and symplectic numerical integrator. The selected Simpson quadrature \eqref{eq:simpsonQuadrature} corresponds to a Lobatto-type quadrature that uses three interpolation nodes (see \eg\ \cite{jay2015}). Simpson's quadrature has been previously used to build variational integrators (\eg, \cite{farr2007,ober2014,rojas2024axioms,capobianco2024}). The work of \cite{farr2007} is centered around a predictive algorithm with adaptive time-steps that conserves angular momentum at fifth order; the fourth-order integrator is applied to the simulation of the $N$-body problem. Reference \cite{rojas2024axioms} precedes this document; it focuses on linear systems and proposes a compact formulation of the variational integrator using partitioned matrices. References \cite{ober2014,capobianco2024} analyze the larger family of Lobatto variational integrators in which the Simpson quadrature leads to a fourth-order method. Examples with separable systems are given in \cite{ober2014}, and \cite{capobianco2024} proposes an event-capturing method, focusing on discontinuous systems involving frictional contact.

Simpson's quadrature can be regarded as part of the larger families of Lobatto IIIA and Lobatto IIIB collocation methods, in their respective fourth-order variants, which are symmetric integrators. Taken individually, these fail to preserve quadratic invariants and are not symplectic, %(see Theorems 2.2 (ch.\,II) and 4.3 (ch.\,VI) in \cite{hairer2006}). However, 
the partitioned Runge-Kutta method based on the Lobatto IIIA--IIIB pairs is symplectic %(see Theorem 4.5 of 
\cite{hairer2006}. This Runge-Kutta method combines one scheme from each of the IIIA and IIIB families to solve separable systems (hence the partitioning of the Hamiltonian) as an implicit method \cite{jay1996}. %Ultimately, the proposed variational integrator from this document and the Runge-Kutta method based on the Lobatto IIIA-IIIB pairs are equivalent as indicated in \cite{marsden2001}. 
We remain in the framework of classical variational integrators \cite{marsden2001,leok2005}. Runge-Kutta methods are not in the scope of this work. In this work, it is the discretization of the action by the discrete Lagrangian that ensures the symplectic property of the integrator \cite{marsden2001,hairer2006}, and it has been remarked that Lobatto quadratures %are well adapted to tackle systems affected by stiffness and these quadratures 
preserve the symmetry of the discrete Lagrangian \cite{marsden2001}. %(see Theorem 6.1 of ch.\,VI in \cite{hairer2006}). 
%The proposed formulation might appeal to those who are familiar with classical mechanics but not with Runge-Kutta methods.

The Simpson variational integrator for nonlinear systems is presented in the next section. After introducing the internal interpolation scheme, the action is discretized to obtain the discrete Euler-Lagrange equations. The complete Simpson scheme, along with the Jacobian matrix involved in Newton's method, is provided in appendix \ref{app:simpson}. Additionally, a nonlinear formulation of the implicit midpoint method \cite{simo1992,marsden2001} is given in appendix \ref{app:Newmark} for reference.  Section \ref{sec:NLDP} presents results for the chaotic, nonlinear double pendulum, regarded as a multibody system. Section \ref{sec:HST} presents numerical results on the Lagrange top to illustrate the accuracy and performance of the proposed integrator. The exact but non-trivial solution is recalled for reference, and convergence results are given with respect to this analytical solution for a motion with loops (section \ref{subsec:completeModel}). The convergence order is further numerically analyzed using a reduced Lagrange top system obtained by selecting special initial conditions for a cuspidal motion. This reduced system is solved using a classical but non-symplectic fourth-order method. In turn, the Simpson variational integrator solves the complete system. The convergence of the difference between both solutions is analyzed and coincides with the convergence order of the more classical method (section \ref{subsec:reducedModel}). Conclusions and perspectives are then discussed in section \ref{sec:conclusion}.

\section{Simpson's variational integrator for nonlinear systems}
%
%
%
%
%
%
%
%
%****************************************************
\subsection{Nonlinear mechanical systems}

Let us describe the targeted systems. Consider a mechanical system with configuration manifold $Q$ denoting the set of states $q\in \mathbb{R}^n$. The velocity phase space is $TQ$, and the Lagrangian is a map $L: TQ \rightarrow \mathbb{R}$. The Lagrangian system is defined by
\begin{equation}
	\label{eq:contLagrangian}
	\mathcal{L}(q,\dot q) = \frac12 \dot q\tran M(q) \dot q - V(q),
\end{equation}
where the overdot implies time differentiation; $M(q)$ is a symmetric positive-definite bilinear form and $V(q)$ is a potential function.

Lagrangian mechanics typically involves the principle of least action \cite{arnold1989,goldstein2002}. Considering the continuous action
\begin{equation*}
	\mathcal{S} = \int_0^T \mathcal L(\dot q, q)\dl t,
\end{equation*}
the motion of the Lagrangian system is such that variations of $\mathcal S$ are zero for an arbitrary variation of the curve $q(t)$, while holding the endpoints of $q(t)$ fixed. The variational procedure results in the well-known Euler-Lagrange equations \cite{arnold1989,goldstein2002} constraining the path followed by the system:
\begin{equation*}
	\diff*{\diffp{L}{\dot q}[]}{t} - \diffp Lq = 0.
\end{equation*}
Taking the Lagrangian \eqref{eq:contLagrangian}, the traditional motion equation for multi-body systems is obtained as
\begin{equation}
	\label{eq:multibodyDynamics}
	M(q)\ddot q + \dot M(q) \dot q - \frac12 \dot q\tran \nabla M(q) \dot q + \nabla V(q) = 0,
\end{equation}
where the second and third terms describe Coriolis and centrifugal effects \cite{spong2020}. More details on the derivatives of $M(q)$ and their computation can be found in \cite{garofalo2013,muller2021}. The proposed integrator is formulated specifically to integrate systems described by equation \eqref{eq:multibodyDynamics}.

%
%
%
%
%
%
%
%
%****************************************************
\subsection{Simpson quadrature and polynomial discretization}

The Simpson discretization uses a quadratic finite elements \cite{raviart1983,allaire2007} internal interpolation at each time step for $0\leqslant t \leqslant h$, with basis functions defined for $0\leqslant \theta \leqslant 1$:
\begin{equation}
	\label{eq:quadraticFE}
	\varphi_0(\theta) = (1-\theta)(1-2\theta),\ 
	\varphi_{\sfrac{1}{2}} (\theta) = 4\theta (1-\theta),\ 
	\varphi_1(\theta) = \theta (2\theta - 1).
\end{equation}
States $q^\alpha(t)\in \mathbb{P}_2$ are approximated with these functions as
\begin{equation}
	\label{eq:approximatedState}
	q^\alpha(t) = \varphi_0 (\theta)q^\alpha_\ell  + \varphi_{\sfrac 12}(\theta)q^\alpha_m  + \varphi_1(\theta)q^\alpha\, ,
\end{equation}
where the subscripts $(\ell,m)$ denote discrete ``left'' and ``middle''values on a specific interval, and the state $q\in Q\subseteq\mathbb{R}^n$ denotes the ``right'' value on each interval. States $q_\ell$ and $q$ are separated by a time step $h\in\mathbb{R}$. On $[0,h]$, $q_\ell = q(0)$, $q_m = q\left(\sfrac h2 \right)$ and $q= q(h)$. Upper Greek indices take as many values as there are configuration parameters: $q\equiv \{ q^\alpha\}$ for $\alpha=1,2,\ldots,n$. By differentating \eqref{eq:quadraticFE},
\begin{equation*}
	\begin{aligned}
		\diff{q}t & = \frac 1h \left( 
			\diff{\varphi_0}{\theta}q_\ell 
			+ \diff{\varphi_{\sfrac 12}}{\theta} q_m 
			+ \diff{\varphi_1}{\theta}q 
			\right)
%		\\
%		& 
			= \frac 1h \left[ 
			(4\theta-3) q_\ell 
			+ 4 (1-2\theta)q_m + (4\theta-1)q 
		\right]
		\\
		& = (1-\theta)g_\ell  + \theta g\, ,
	\end{aligned}
\end{equation*}
where discrete velocities $g_\ell, g_m, g \in Q\times Q$ are given by Gear's scheme \cite{gear1971}:
\begin{equation*}
	g_\ell^\alpha 
	%& 
	= \frac 1h \left(-3q_\ell^\alpha + 4 q^\alpha_m - q^\alpha \right),
	%\\
	\quad
	g_m^\alpha 
	%& 
	= \frac1h \left(-q_\ell^\alpha + q^\alpha \right),
	%\\
	\quad
	g^\alpha
	%&
	= \frac1h \left(q_\ell^\alpha - 4q_m^\alpha + 3q^\alpha \right).
\end{equation*}

The discrete Lagrangian approximates the action integral $\mathcal S$ along the curve segment between $q_\ell$ and $q$. The approximation to $\int_0^T \mathcal L\dl t$ is given by the Simpson quadrature \eqref{eq:simpsonQuadrature}, as in \cite{dubois2023_GSI,dubois2023_FVCA,rojas2024axioms}:
\begin{equation}
	\int_0^h \mathcal{L}(q(\theta),\dot q(\theta))\dl \theta \simeq \frac h6 \left(\mathcal L(0) + 4\mathcal L \left(\frac 12\right) + \mathcal L(1)\right).
\label{eq:simpsonQuadrature}
\end{equation}
As prescribed in \cite{leok2005}, the discrete action becomes the fundamental object. It is the sum of discrete Lagrangians, $L_\mathrm{d}:Q\times Q \rightarrow \mathbb{R}$, approximates the exact discrete Lagrangian and is a generating function of the symplectic flow as remarked by \cite{marsden2001,leok2005,hairer2006}:
\begin{multline}
	L_\mathrm{d} 
	= \frac h2\left[ 
		\frac16 M_{\alpha\beta}( q_\ell )\, g_\ell^\alpha g_\ell^\beta 
		+ \frac23 M_{\alpha\beta}( q_m )\, g_m^\alpha g_m^\beta
		+ \frac16 M_{\alpha\beta}(q)\, g^\alpha g^\beta
\right]
	\\
	- h\left[ 
		\frac16 V( q_\ell )
		+\frac23 V( q_m )
		+\frac16 V(q)
	\right]\, .
\end{multline}
% . 
Considering a discrete curve of points $\{q_0,\ldots,q_N\}$, the discrete action, for a motion $t\mapsto q(t)$ with $0\leq t\leq T$, is a sum $S_\mathrm{d}:Q^{n+1}\rightarrow \mathbb{R}$ of discrete Lagrangians over a finite number $N$ of discretization intervals:
\begin{equation}
\label{eq:actionSum}
	S_\mathrm{d} = \ldots + L_\mathrm{d}(q_{j-1},q_{j-\sfrac12},q_j) +  L_\mathrm{d}(q_j,q_{j+\sfrac12},q_{j+1}) + \ldots.
\end{equation}
Applying the discrete variational principle $\delta S_\mathrm{d} = 0$, for arbitrary variations of the states, results in the discrete Euler-Lagrange equations. At the middle of the interval, $\delta S_\mathrm{d} = 0$ for an arbitrary variation $\delta q_{j+\sfrac12}$ implies that 
\begin{equation}
	\label{eq:SimpELmiddle}
	\diffp{L_\mathrm{d}}{q_m^\gamma} 
	= \frac13 M_{\gamma\beta}( q_\ell )\,g_\ell^\beta
		-\frac13 M_{\gamma\beta}(q)\,g^\beta
		+\frac h6 \partial_\gamma M_{\alpha\beta}( q_m )\,g_m^\alpha g_m^\beta
		-\frac h3 \partial_\gamma V( q_m )
	= 0\, .
\end{equation}
On the borders of each interval, the discrete Euler-Lagrange equations are obtained when $\delta S_\mathrm{d} = 0$ for an arbitrary variation $\delta q_j$, that is,
\begin{equation}
	\diffp{L_d}{q^\gamma}(q_{j-1},q_{j-\sfrac12},q_j) + \diffp{L_d}{q_\ell^\gamma}(q_j,q_{j+\sfrac12},q_{j+1}) = 0\, 
	\label{eq:DEL}
\end{equation}
Following this procedure, the resulting variational integrator is guaranteed to be symplectic \cite{marsden2001,hairer2006}.

%%%%%%
%%%%%%
%%%%%%
\subsection{Simpson variational integrator}

%*** Voir manuscrit de François, datant du 18 janvier

The terms of the discrete Euler-Lagrange equations \eqref{eq:DEL} can be identified as being a sum of discrete generalized momenta. Using the shorthand $\partial_\alpha \equiv \diffp{}{q^\alpha}$, the momenta, on the right of the interval, are defined by
\begin{equation}
	\label{eq:momentaRight}
	\begin{aligned}
	p_{\gamma}  
	& = \diffp{L_\mathrm{d}}{q^\gamma}(q_\ell,q_m,q)
	\\
	& = -\frac16 M_{\gamma\beta}( q_\ell )\,g_\ell^\beta
		+ \frac23 M_{\gamma\beta}( q_m )\,g_m^\beta
		+ \frac12 M_{\gamma\beta}(q)\,g^\beta
		+\frac{h}{12}\partial_\gamma M_{\alpha\beta}(q)\,g^\alpha g^\beta 
		- \frac h6 \partial_\gamma V(q).
	\end{aligned}
\end{equation}
Adding half of the right-hand side of equation \eqref{eq:SimpELmiddle}, to the momenta expressed above, yields:
\begin{multline}
	\label{eq:SimpRightMomenta}
	p_{\gamma} = \frac16 M_{\gamma\beta}( q_\ell )\,g_\ell^\beta
		+ \frac23 M_{\gamma\beta}( q_m )\,g_m^\beta
		+ \frac16 M_{\gamma\beta}(q)\,g^\beta
	\\
		+ \frac{h}{6}\partial_\gamma M_{\alpha\beta}( q_m )\,g_m^\alpha g_m^\beta
		+ \frac{h}{12}\partial_\gamma M_{\alpha\beta}(q)\,g^\alpha g^\beta 
		- \frac h3 \partial_\gamma V( q_m )\,
		- \frac h6 \partial_\gamma V(q)\, .
\end{multline}
Using the definition \eqref{eq:SimpRightMomenta}, the discrete Euler-Lagrange equations on the border of the interval become
\begin{equation}
	\label{eq:discreteELeq}
	\diffp{L_d}{q^\gamma}(q_{j-1},q_{j-\sfrac12},q_j) + \diffp{L_d}{q_\ell^\gamma}(q_j,q_{j+\sfrac12},q_{j+1}) = 0\,.
\end{equation}
The term on the left is, by definition, $p_{\gamma,j}$. By \eqref{eq:discreteELeq},  
\begin{equation}
	\label{eq:jMomenta}
	p_{\gamma,j} = - \diffp{L_d}{q_\ell^\gamma}(q_j,q_{j+\sfrac12},q_{j+1})\,.
\end{equation}
Also, by \eqref{eq:momentaRight},
\begin{equation}
	\label{eq:jp1Momenta}
	p_{\gamma,j+1} = \diffp{L_d}{q_r^\gamma}(q_j,q_{j+\sfrac12},q_{j+1})\,.
\end{equation}

In the above, the discrete Lagrangian involves an internal configuration. For the linear cases, this internal configuration can be eliminated, leading to a reduced discrete Lagrangian. This reduced Lagrangian has been established explicitly in \cite{rojas2024axioms}. In the present nonlinear case, an explicit expression of this reduced Lagrangian cannot be written, but exists by the implicit function theorem. In this logic and as summarized by \cite{hairer2006}, \eqref{eq:jMomenta} defines a bijection between $p_j$ and $q_{j+1}$ for a given $q_j$. The resulting method maps $(p_j,q_j) \mapsto (p_{j+1},q_{j+1})$ and is symplectic (see %Theorems 5.1 and 6.1 of 
ch.\,VI in \cite{hairer2006}, or part two of \cite{marsden2001}). The following sections apply the proposed integrator to examples of nonlinear systems that at least preserve one physical quantity. 

Subtracting equations \eqref{eq:jMomenta} and \eqref{eq:jp1Momenta} gives the discrete evolution of momenta. Analogously, adding \eqref{eq:jMomenta} and \eqref{eq:jp1Momenta} gives the discrete evolution of the configurations. Combining these two discrete equations with the discrete Euler-Lagrange equation on the internal configuration \eqref{eq:SimpELmiddle} leads to the nonlinear Simpson scheme for Lagrangian systems. Detailed expressions are provided in Appendix \ref{app:Newmark}, along with the Jacobian matrix required to initialize Newton's method for solving the nonlinear problem.

\section{Nonlinear double pendulum}
\label{sec:NLDP}

This section presents some results on the nonlinear double pendulum system. This multibody system is chaotic. Only the system energy is conserved.

\subsection{System description}
\label{subsec:dpDescription}

\begin{figure}[H]
	\centering
	\includegraphics{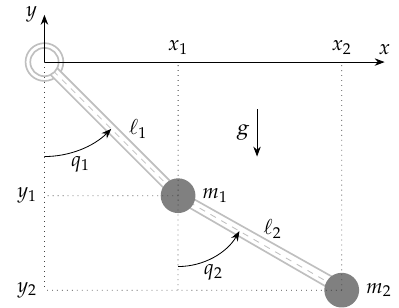}
	\caption{Double pendulum affected by gravity. Two point masses $\{m_1, m_2\}$ are linked together by massless thin rigid rods of respective lengths $\{\ell_1,\ell_2\}$. Mass positions are given by the generalized coordinates $q=\left(q^1,q^2\right)$.}
	\label{fig:doublePendulum}
\end{figure}

The double pendulum shown in Figure \ref{fig:doublePendulum} is a two-degree-of-freedom nonlinear system. It is formed by two point masses $\{m_1, m_2\}$ linked together by massless thin rigid rods of respective lengths $\{\ell_1,\ell_2\}$. Each joint brings one degree of freedom to the articulated system. 

The masses' locations are given by the generalized coordinates $q=\left(q^1,q^2\right)$ as
\begin{equation*}
\begin{aligned}
	\left(x^1,y^1\right) & = \left(\ell_1 \sin q^1,-\ell_1 \cos q^1\right)\\
	\left(x^2,y^2\right) & = \left(\ell_1 \sin q^1 + \ell_2 \sin q^2,-\ell_1 \cos q^1 - \ell_2 \cos q^2\right).
\end{aligned}
\end{equation*}
Velocities are obtained by time differentiation, considering that $q_i\equiv q_i(t)$. The kinetic energy  is 
\begin{equation*}
		T 
		= \frac12 (m_1+m_2){\ell_1}^2 \left(\dot{q}^1\right)^2 + \frac 12 m_2 {\ell_2}^2 \left(\dot{q}^2\right)^2 + m_2 \ell_1 \ell_2 \dot{q}^1\dot{q}^2\cos\left(q^1-q^2\right),
\end{equation*}
where an overdot indicates time differentiation. The kinetic energy being defined as $T = \frac12 \dot q\tran M \dot q$, the mass matrix $M(q)$ (symmetric and positive-definite) is identified from the above as
\begin{equation*}
	M = \begin{pmatrix}
		(m_1+m_2){\ell_1}^2
		& m_2 \ell_1 \ell_2 \cos(q_1-q_2)
		\\[1ex]
		m_2 \ell_1 \ell_2 \cos(q_1-q_2)
		& m_2 {\ell_2}^2
	\end{pmatrix}.
\end{equation*}
Taking $\cg$ as the gravity constant, the potential energy is
\begin{equation*}
	V= - m_1 \cg  \ell_1 \cos q_1 - m_2 \cg \left(\ell_1\cos q_1 + \ell_2\cos q_2\right)\,.
\end{equation*}

By regarding this double pendulum as a conservative system, the Hamiltonian function is established as the total energy of the system: 
\begin{equation}
	\label{eq:Hamiltonian}
	H(p,q) 
	= 
	\frac12 p^T M(q)^{-1} p + V(q)\,.
\end{equation}
Finally, motion is described by the canonical equations $\dot p_i = -\diffp{H}{q^i}$ and $\dot q^i = \diffp{H}{p_i}$.

%
%
%
%***********************
\subsection{Numerical experiments}

The double pendulum system described in this section is chaotic. Only the system energy is a conserved quantity; the system is not integrable, and its Hamiltonian is inseparable. Because of these characteristics, no exact solutions can be used as a reference. In the following analysis, the classical fourth-order Runge-Kutta method (RK4, see \eg\ \cite{suli2003}) is used to obtain reference curves, with a small enough step size to consider results as reasonably correct. For further reference, results obtained with the symplectic Implicit midpoint method \cite{simo1992} are also provided. All computations in this section were conducted with arbitrary precision and 64 digits, using the \texttt{mpmath} Python library \cite{mpmath2023}. 

Table \ref{tab:DPsimulationParams} shows the constants and initial conditions used for the simulations in the current section. The link lengths are chosen so that $\ell_1 = \ell_2 \equiv \ell = \sfrac{\cg}{\omega_0}$, where $\cg=\SI{9.81}{\meter\per\second\squared}$ is the acceleration of the gravity field. 
\begin{table}[H] 
\caption{
	Constants and initial conditions used for numerical simulations.
	\label{tab:DPsimulationParams}
	}
\newcolumntype{C}{>{\centering\arraybackslash}X}
\begin{tabularx}{\textwidth}{CCCC}
	\toprule
	\multicolumn{2}{c}{\textbf{Constants}}	
	& \multicolumn{2}{c}{\textbf{Initial conditions}}
	\\
	\cmidrule(lr){1-2}
	\cmidrule(lr){3-4}
	$m_1=m_2$
	& \SI{1}{\kilogram}		
	& ${q}(0)$ 
	& $(\sfrac \pi4,\sfrac \pi3)^T \,\si{\radian}$
	\\[0.5ex]
	$\omega_0$ 
	& $2\pi\,\si{\per\second}$		
	& ${p}(0)$ 
	& $(0,0)^T \, \si{\kilogram\meter\squared\per\second}$
	\\
\bottomrule
\end{tabularx}
\end{table}
\subsubsection{States evolution}

\begin{figure}[H]
	\centering
	\begin{adjustwidth}{-\extralength}{-\extralength}
	\includegraphics{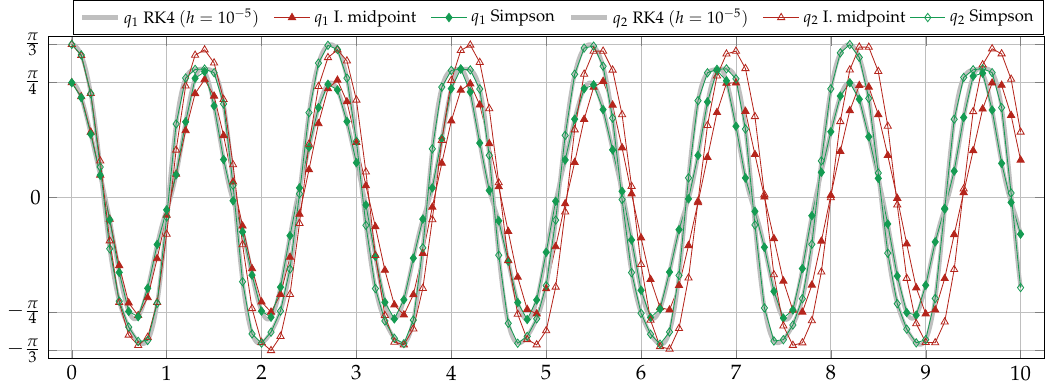}
	\caption{Configuration parameters ($q$) evolution for the nonlinear double pendulum. Step size is $h=\num{e-1}$ for the symplectic integrators (Implicit midpoint and Simpson). RK4 gives the reference curve at $h=\num{e-5}$. Simpson's solutions follow the reference curve for longer simulations.}
	\end{adjustwidth}
	\label{fig:statesDP}
\end{figure}

\begin{figure}[t!]
	\centering
	\begin{adjustwidth}{-\extralength}{-\extralength}
	\includegraphics{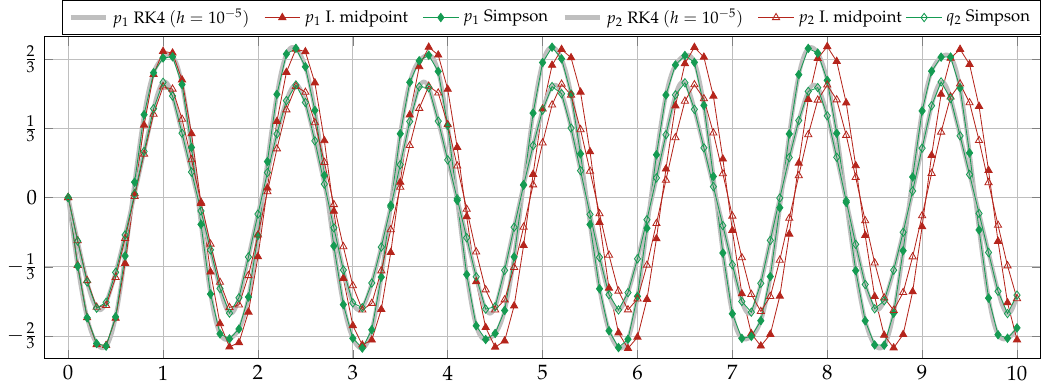}
	\caption{Generalized momenta ($p$) evolution for the nonlinear double pendulum. %Initial conditions are specified in Table \ref{tab:DPsimulationParams}. 
Step size is $h=\num{e-1}$ for the symplectic integrators (Implicit midpoint and Simpson). RK4 gives the reference curve at $h=\num{e-5}$. Simpson's solutions follow the reference curve for longer simulations.}
	\end{adjustwidth}
	\label{fig:momentaDP}
\end{figure}

Figures \ref{fig:statesDP} and \ref{fig:momentaDP} show the configuration parameters and momenta evolution, respectively, for the nonlinear double pendulum motion approximation by three methods. Since the RK4 method is widely used and known, the reference curve (in light gray color) was obtained using RK4 at a small step size $h=\num{e-5}$, and the other curves (at $h=\num{e-1}$) are compared with this one. Early on, the Implicit midpoint and RK4 (at $h=\num{e-1}$) deviate from the reference curve. However, it can be appreciated that the Simpson curves correctly follow the reference curve across the simulation time window of \SI{10}{\second}.

\subsubsection{Energy convergence}

\begin{figure}[H]
	\begin{adjustwidth}{-\extralength}{-\extralength}
	\centering
		\includegraphics{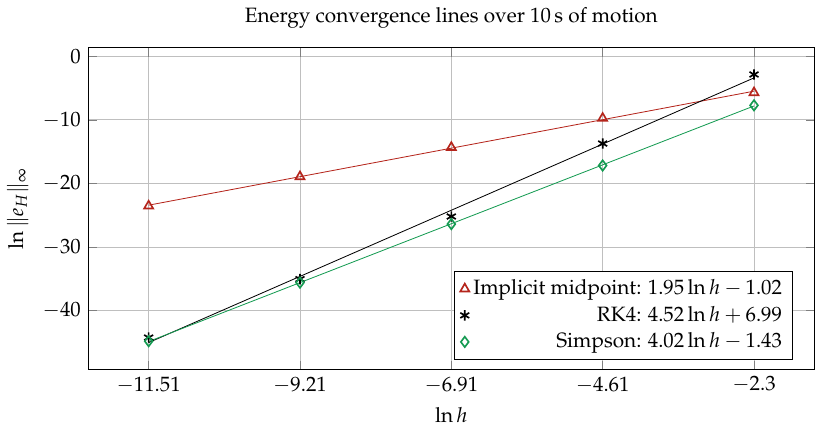}
	\end{adjustwidth}
	\caption{Energy error convergence lines for the double pendulum motion. RK4 showcases a higher convergence rate than its characteristic one (of 4) because the error is high for the larger time steps. The convergence rates of the Implicit midpoint and Simpson methods are 2 and 4, respectively.}
	\label{fig:DPenergyConvergence}
\end{figure}

The energy quantity $H(p,q)$ is the only motion constant for this system. A numerical analysis on convergence towards this quantity illustrates the performance of the proposed integrator. The energy error is calculated according to $e_H= \frac{H(p,q) -H_0}{H_0}$, where $H_0 = H(p_0,q_0)$ with the initial conditions of Table \ref{tab:DPsimulationParams}. The convergence rate is measured according to the procedure found in \cite{allaire2007}; precision is evaluated on the $\ell^\infty$ error norm
$\|e_H\|_\infty = \sup_n \left\lvert\frac{H(p_n,q_n)-H_0}{H_0}\right\lvert$. Errors are calculated for decreasing values of the step size, computing the $\|e_H\|_\infty$ norm for each case. These errors are plotted in Figure \ref{fig:DPenergyConvergence}, on a logarithmic scale.

\begin{figure}[H]
	\begin{adjustwidth}{-\extralength}{-\extralength}
	\centering
		\includegraphics{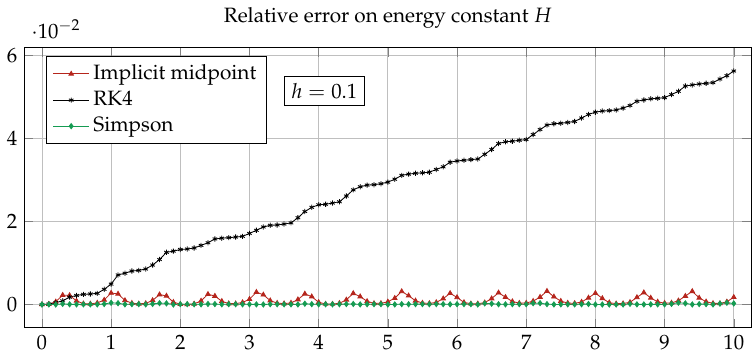}
	\end{adjustwidth}
	\caption{Relative error on the conserved quantity $H(p,q)$ for the double pendulum. RK4's error grows with time. The symplectic methods (Implicit midpoint and Simpson) do not artificially dissipate energy.}
	\label{fig:energyErrorEvolution}
\end{figure}

The Implicit midpoint shows its characteristic second-order convergence rate. The RK4 method shows an increased convergence rate due to the fact that the error values obtained for larger time steps are high. The Simpson method has a fourth-order convergence rate. Figure \ref{fig:energyErrorEvolution} shows the evolution of the error on the energy constant $H$ across \SI{10}{\second} of motion. As expected from the variational integrators, no artificial energy dissipation occurs, with the Simpson approximation remaining closest to zero. Table \ref{tab:convergenceEvolution} shows the evolution of the energy error norm relative to the increase in simulation time. The symplectic integrators (Implicit midpoint and Simpson) hold their respective convergence rates across all tests, regardless of the increased simulation time. Conversely, the convergence order of the RK4 method suffers from fluctuations across tests. 

\begin{table}[H]
	\sisetup{scientific-notation = true,round-mode=places,round-precision=2}
	\caption{Convergence order with respect to simulation length for motion simulations held on a nonlinear double pendulum. Initial conditions are specified in Table \ref{tab:DPsimulationParams}. The error norm $\|e_H\|_\infty$ grows over time on the non-symplectic method RK4. The symplectic methods (Implicit midpoint and Simpson) preserve their convergence order on energy.
	}
	\label{tab:convergenceEvolution}
	\begin{adjustwidth}{-\extralength}{-\extralength}
	\centering
		\includegraphics[]{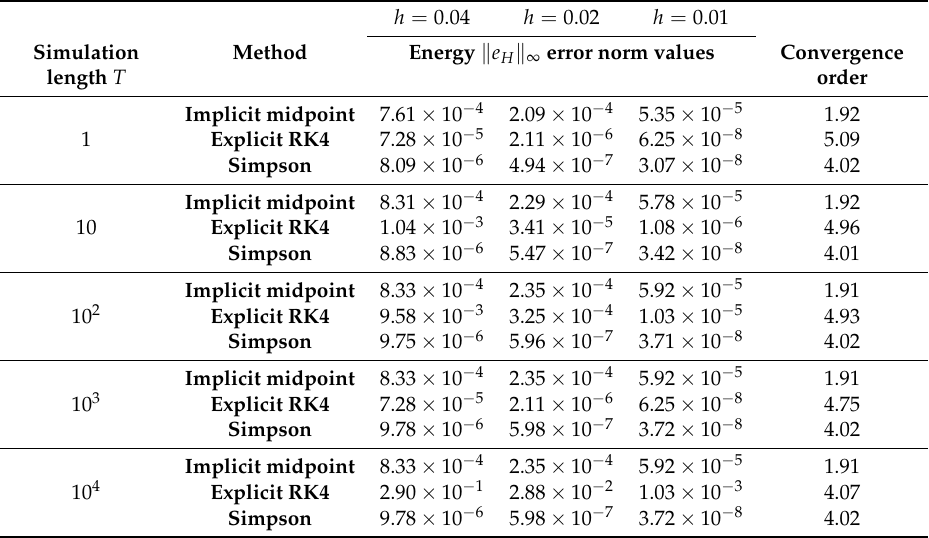}
	\end{adjustwidth}
\end{table}

%
%
%
%
%
%
%
%
%****************************************************
\section{The Lagrange top}
\label{sec:HST}

We now apply the Simpson variational integrator to the Lagrange top. It is a three-degrees-of-freedom nonlinear system for which the analytical solution is nontrivial. It is a good example to illustrate the performance of the proposed integrator.

\subsection{System description}
\label{subsec:topDescription}

The Lagrange top is a heavy symmetric top of mass $m_\mathrm{top}$, which is affected by gravity, and spins on a table with a fixed constant point. It is a rigid body with cylindrical symmetry around one of its axes, so that two of its principal moments of inertia are equal. Let us choose the body symmetry axis to be $z$. Its moments of inertia are $I_1=I_2\equiv I$, and $I_3$, and the point of contact between the top and the table is fixed and allows for rotational motion in an idealized manner. 

\begin{figure}[H]
	\centering
	\includegraphics[scale=1]{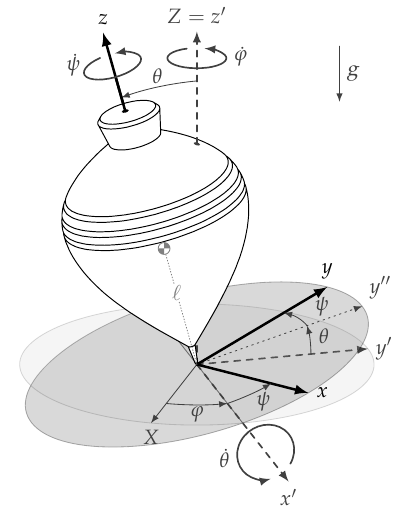}
	\caption{Euler angles ``ZXZ'' sequence for the Lagrange top. Three angles $q=(\varphi,\theta,\psi)$ characterize the spinning motion. In order, $(\psi,\theta,\varphi)$ respectively describe: the top spin about its symmetry axis $z$; the $z$ body axis precession about the inertial frame axis $Z$; the nutation of the $z$ body axis about the inertial frame axis $Z$. The center of mass is located along $z$ at a distance $\ell$ from the contact point. Gravity points downwards along the vertical direction.}
	\label{fig:HST}
\end{figure}

The body motion equations are derived as seen from an inertial frame of reference, external to the body, following the Euler-Lagrange equations. Let us use Euler angles $q=\left(\varphi,\theta,\psi\right)$ as the spinning top motion coordinates. The specific sequence is shown in Figure \ref{fig:HST}. Using the Euler angles description of Figure \ref{fig:HST} results in angular velocities taking place across different frames.  Spinning takes place on the body  frame while precession takes place on the fixed reference frame, yielding the body's angular velocity
\begin{equation*}
	\omega = \dot\psi\, z + \dot\theta\, x^\prime + \dot\varphi\, Z.
\end{equation*}
By projecting $\omega$ entirely on the body basis vectors $(x,y,z)$, its expression becomes
\begin{equation*}
	\omega = \dot\theta\,x + \dot\varphi\sin\theta\,y  + \left(\dot\psi + \dot\varphi\cos\theta\right) z,
\end{equation*}
because $Z = \cos\theta\,z + \sin\theta y$ and nutation occurs around $x$ during motion. Note that we are assuming a symmetrical top with $I_1=I_2\equiv I$. Therefore, any axes in the $(x,y)$ plane are principal axes (such as the line of nodes $x^\prime$).

The kinetic energy of the top is given by
\begin{equation}
	\label{eq:KEtop}
	T = \frac12 \omega\tran J \omega = \frac12 I_3\left(\dot\psi + \dot\varphi\cos\theta\right)^2 + \frac12 I\left(\dot\varphi^2\sin^2\theta + \dot\theta^2\right),
\end{equation}
where $J = \diag (I,I,I_3)$. The system mass tensor can be identified from the above as
\begin{equation*}
	M(q) = \begin{pmatrix}
		I\sin^2\theta + I_3 \cos^2\theta
		& 0 
		& I_3\cos\theta
		\\
		0
		& I
		& 0
		\\
		I_3\cos\theta
		& 0
		& I_3
	\end{pmatrix}\,,
\end{equation*}
so the kinetic energy can be equivalently written as $T = \frac12 \dot q\tran M(q)\dot q$. 
The potential energy is
\begin{equation}
	\label{eq:PEtop}
	V = m_\mathrm{top}\cg\ell\cos\theta \equiv I m \cos\theta,
\end{equation}
where $\ell$ is the distance from the fixed point of contact between the table and the top to the center of mass of the top. The parameter $m$ is  
given by $m=\sfrac{m_\mathrm{top}\cg\ell}{I}$.

Recall that the Lagrangian is given by $L=T-V$ and $q=\left(\varphi,\theta,\psi\right)$; using the Euler-Lagrange equations, three equations of motion can be deduced:
\begin{subequations}
\label{eq:ELtop}
\begin{align}
	\label{eq:ELtop1}
	\diff{p_\varphi}{t}
	& = 0
	\\
	\label{eq:ELtop2}
	I\ddot\theta - \left( I m + I \dot\varphi^2\cos\theta - I_3\left(\dot\psi+\dot\varphi\cos\theta\right)\dot\varphi \right)\sin\theta
	& = 0,	
	\\
	\label{eq:ELtop3}
	\diff{p_\psi}{t}
	& = 0,
\end{align}
\end{subequations}
where $\left(p_\varphi,p_\psi\right)$ are two conserved generalized momenta defined by
\begin{subequations}
\label{eq:conservedMomenta}
\begin{align}
	\label{eq:pPhi}
	p_\varphi = \diffp{L}{\dot\varphi}
	& = I_3\left(\dot\psi+\dot\varphi\cos\theta\right)\cos\theta + I\dot\varphi\sin^2\theta,	
	\\
	\label{eq:pPsi}
	p_\psi = \diffp{L}{\dot\psi}
	& = I_3\left(\dot\psi+\dot\varphi\cos\theta\right).
\end{align}
\end{subequations}
Equations \eqref{eq:ELtop1} and \eqref{eq:ELtop3} describe two conserved quantities relating to the cyclic coordinates $(\varphi,\psi)$. These are the angular momenta in the $Z$ and $z$ directions, respectively. More fundamentally, these momenta are conserved because no torque is being exerted in the $(Z,z)$ plane, but only along the $x$ direction. The system \eqref{eq:ELtop}--\eqref{eq:conservedMomenta} shall be called \emph{complete Lagrange top} system afterwards.

\subsection{Comparison with the exact solution of the complete Lagrange top system}
\label{subsec:completeModel}

The Simpson variational integrator is now evaluated on the Lagrange top motion. For further reference, results obtained with other widely used numerical methods are also provided: the Implicit midpoint method, which is symplectic \cite{simo1992}, and the explicit Runge-Kutta's fourth-order method (RK4), which is a standard method. All computations in this section were conducted with arbitrary precision and 64 digits, using the \texttt{mpmath} Python library \cite{mpmath2023}.

The simulation parameters of this section are listed in Table \ref{tab:HSTparameters}. The physical parameters of the top are taken from \cite{laos2022} and come from an existing spinning top toy.

\begin{table}[H] 
\sisetup{
	scientific-notation = true,
	round-mode = places,
	round-precision = 2}
\caption{Lagrange top parameters for a motion showcasing loops.\label{tab:HSTparameters}}
\newcolumntype{C}{>{\centering\arraybackslash}X}
\begin{tabularx}{\textwidth}{p{3cm}CCC}
	\toprule
	
	& %\multicolumn{2}{C}
	{\textbf{Parameter}}
	& \textbf{Value}
	& \textbf{Unit}
	\\
	\midrule
	\multirow{4}{=}{Physical parameters of the top \cite{laos2022}}
	& $m_\mathrm{top}$
	& \num{0.1}
	& [\si{\kilogram}]
	\\
%	Cross-sectional moment of inertia 
	& $I$
	& \num{0.002329969592394382}
	& [\si{\kilogram\meter\squared}]
	\\
%	Longitudinal moment of inertia
	& $I_3$
	& \num{0.000125}
	& [\si{\kilogram\meter\squared}]
	\\
%	Center of mass location
	& $\ell$
	& \num{0.15}
	& [\si{\meter}]
	\\
	\midrule
	\multirow{3}{*}{Initial conditions}
	& $\left\{\varphi_0,\theta_0,\psi_0\right\}$
	& $\left\{0,\,\sfrac{\pi}{3},\,0\right\}$
	& [\si{\radian}]
	\\[0.5ex]
	& $\left\{\dot\varphi_0,\dot\theta_0,\dot\psi_0\right\}$
	& $\left\{9.2,\,0,\,252\right\}$
	& [\si{\radian\per\second}]
	\\
	\bottomrule
\end{tabularx}
\end{table}

The values in Table \ref{tab:HSTparameters} are applied to equation \eqref{eq:exactTheta} (below), which expresses the exact value of the nutation angle $\theta$ in terms of the time \cite{whittaker1917}: 
\begin{equation}
	\cos\theta 
	= 
	\frac{2}{m} \wp(t+\omega_3) + \frac{2 I c + {p_\psi}^2}{6 I m},
	\label{eq:exactTheta}
\end{equation}
where $m=\sfrac{m_\mathrm{top}\cg\ell}{I}$ according to \eqref{eq:PEtop}; $\wp$ denotes the elliptic Weierstrass $\wp$-function; $\omega_3$ is an imaginary number corresponding to a half-period; $c=\frac12 I {\dot \theta_0}\strut^2+\frac{(p_\phi-p_\psi \cos\theta_0)^2}{2I\sin^2\theta}+Im\cos\theta_0$. Full details on this exact solution can be found in \cite{whittaker1917}. For the studied motion, the period of nutation is $\bar t=\SI{1.84671}{\second}$. We will show our results in terms of this nutation period. Therefore, the effective step size will be expressed as a portion of this period $\bar t$.

\begin{figure}[H]
	\begin{adjustwidth}{-\extralength}{-\extralength}
	\centering
		\includegraphics[]{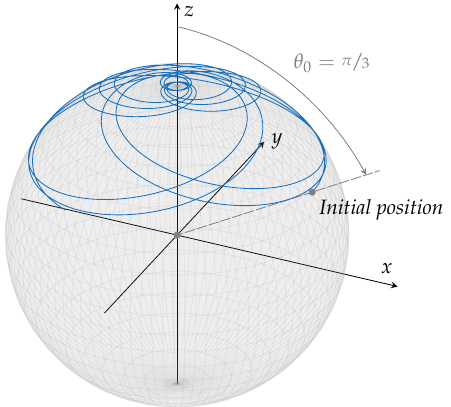}
		\qquad
		\includegraphics[]{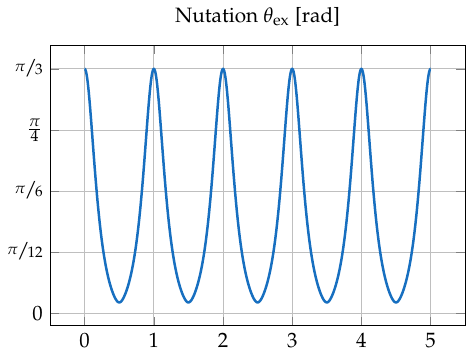}
	\end{adjustwidth}
	\caption{Selected reference top motion over 5 periods of nutation (see  Table \ref{tab:HSTparameters}). This example belongs to the class of ``looping'' motions of the Lagrange top. The initial nutation $\theta_0=\sfrac{\pi}{3}$ bounds motion: $0<\theta\leqslant\theta_0$.}
\end{figure}

\subsubsection{Nutation, energy and momentum}

\begin{figure}[H]
	\begin{adjustwidth}{-\extralength}{-\extralength}
	\centering
		\includegraphics[]{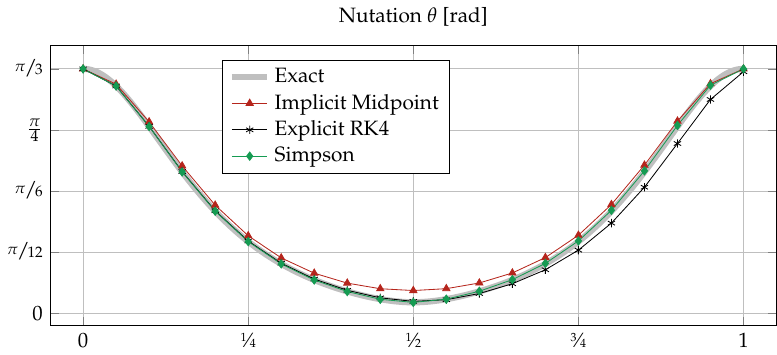}
	\end{adjustwidth}
	\caption{Several methods compared against the exact nutation solution over one period for $h=0.05$. Simpson's solution correctly follows the exact one.}
	\label{fig:nutation}
\end{figure}

Figure \ref{fig:nutation} shows the nutation approximation of three integrators, compared against the exact solution for one period of motion. The Simpson approximation correctly follows the exact solution with no visible deviations.

\begin{figure}[H]
	\begin{adjustwidth}{-\extralength}{-\extralength}
	\centering
		\includegraphics[]{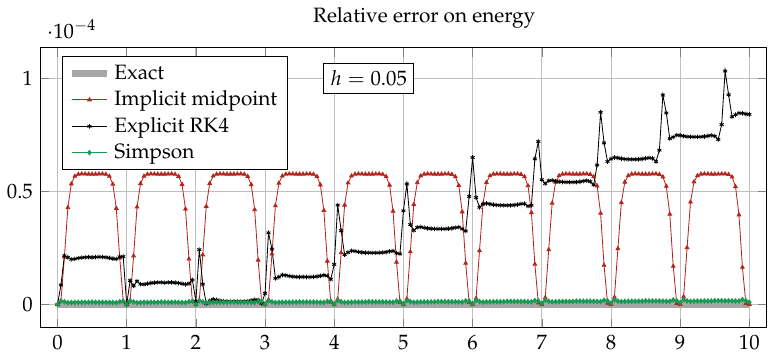}
	\end{adjustwidth}
	\caption{Relative error on the energy constant $H$ across 10 periods of motion. The variational integrators (Implicit midpoint and Simpson) showcase a good energy behavior for reasonably large step sizes and do not artificially dissipate energy (unlike the non-symplectic RK4 method).}
	\label{fig:energyTop}
\end{figure}

Figure \ref{fig:energyTop} shows the evolution of the error on the energy constant $H$ across ten periods of motion. As expected from the variational integrators, the energy is not being artificially dissipated across time, with the Simpson approximation being the closest to zero.

\begin{figure}[H]
	\begin{adjustwidth}{-\extralength}{-\extralength}
	\centering
		\includegraphics[]{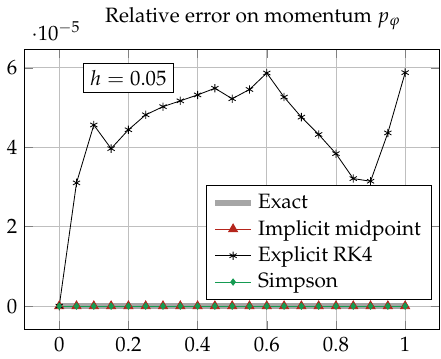}
		\qquad
		\includegraphics[]{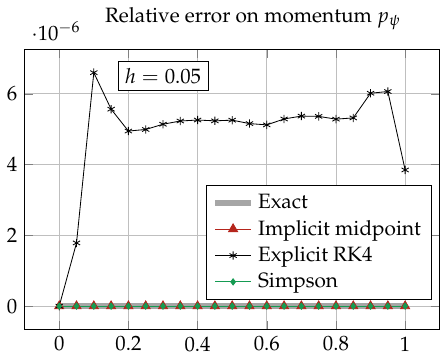}
	\end{adjustwidth}
	\caption{Relative errors on conserved momenta. The symplectic integrators (Implicit midpoint and Simpson) preserve these momenta. 
	}
	\label{fig:MomentumConservation}
\end{figure}

Figure \ref{fig:MomentumConservation} illustrates how the Simpson integrator (as well as the Implicit midpoint method) exactly preserves the conserved momenta $p_\varphi$ and $p_\psi$, as dictated by equations \eqref{eq:ELtop}. 

\subsubsection{Convergence}
\label{subsubsec:convergence}

\begin{figure}[H]
	\begin{adjustwidth}{-\extralength}{-\extralength}
	\centering
		\includegraphics[]{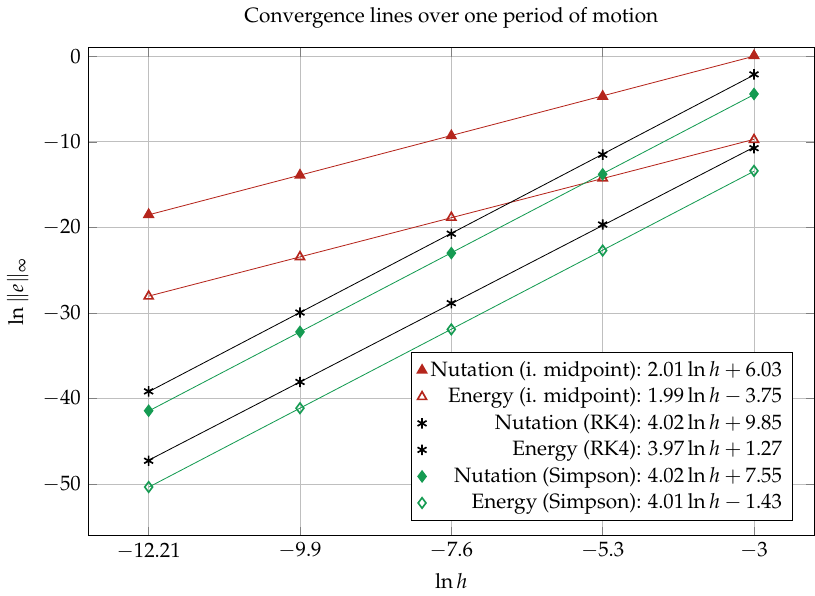}
	\end{adjustwidth}
	\caption{Nutation (red) and energy (green) error convergence rates of various methods. RK4 showcases its characteristic fourth-order convergence, and the implicit midpoint method is second-order convergent. Simpson's variational integrator is fourth-order convergent but more precise than RK4.}
	\label{fig:thetaConvergenceREL}
\end{figure}

The relative error on nutation is calculated according to $	e_\theta(t)= \frac{\theta -\theta_{\text{ex}}}{\theta_{\text{ex}}}$. The energy error is calculated according to $	e_H= \frac{H(p,q) -H_0}{H_0}$, where $H_0 = H(p_0,q_0)$ and is calculated with the initial conditions of Table \ref{tab:HSTparameters}. The convergence rate is measured according to the procedure found in \cite{allaire2007}. The precision of the methods is evaluated on the $\ell^\infty$ error norms
\begin{equation*}
	\|e_\theta\|_\infty 
	= 
	\sup_n \left\lvert \frac{\theta_n-\theta_{\text{ex}_n}}{\theta_{\text{ex}_n}} \right\lvert,
	\qquad
	\|e_H\|_\infty 
	= 
	\sup_n \left\lvert\frac{H(p_n,q_n)-H_0}{H_0}\right\lvert.
\end{equation*}
Errors are calculated for decreasing values of the step size over one period of nutation motion. The $\|e_\theta \|_\infty$ and $\|e_H\|_\infty$ norms were computed for each case. These errors are plotted in Figure \ref{fig:thetaConvergenceREL}, on the logarithmic scale.

Tables \ref{tab:nutationErrorNorms} and \ref{tab:energyErrorNorms} show the evolution of the nutation error norm and energy error norm, respectively, relative to the increase in simulation time.

\begin{table}[H]
	\caption{Nutation error norms and convergence order evolution on increasing simulation lengths of the heavy symmetric top motion (see Figure \ref{fig:HST}). Initial conditions for the experiment are specified in Table \ref{tab:HSTparameters}. Error norms $\|e_\theta\|_\infty$ increase with simulation length $T$. The Simpson integrator convergence order fluctuation is minimal over many simulation periods.\smallskip}
	\label{tab:nutationErrorNorms}
	\begin{adjustwidth}{-\extralength}{-\extralength}
		\centering
		\includegraphics{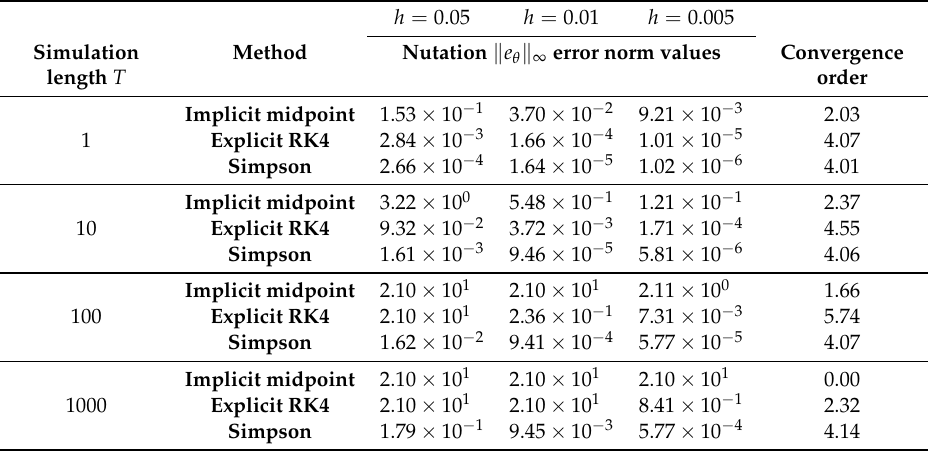}
	\end{adjustwidth}
\end{table}

Regarding nutation (Table \ref{tab:nutationErrorNorms}), the convergence order of the Implicit midpoint method degrades and falls to zero for 1000 periods of simulation because of the accumulated errors over time. The non-symplectic RK4 does not maintain its convergence rate either on the nutation error. This estimated convergence order fluctuates and even approaches the value of 6 for 100 periods of motion, only because the error at $h=0.05$ is very high. Conversely, the Simpson integrator is able to maintain its convergence order for longer simulations by not deviating much from the exact solution.

On the energy error norm evolution (Table \ref{tab:energyErrorNorms}), all methods showcase better results than those obtained for the nutation approximation. The Implicit method convergence rate is stable over time. Fluctuations of the convergence order persist on the RK4 method, with its convergence rate approaching the value of 5, because the error is high for larger time steps. The Simpson method preserves its fourth-order convergence rate for long simulations.

\begin{table}[H]
	\caption{Energy error norms and convergence order evolution on increasing simulation lengths of the heavy symmetric top motion (see Figure \ref{fig:HST}). Initial conditions for the experiment are specified in Table \ref{tab:HSTparameters}. Error norms $\|e_H\|_\infty$ increase with simulation length $T$. The symplectic methods (Implicit midpoint and Simpson) preserve their convergence order on energy.\smallskip}
	\label{tab:energyErrorNorms}
	\begin{adjustwidth}{-\extralength}{-\extralength}
		\centering
		\includegraphics{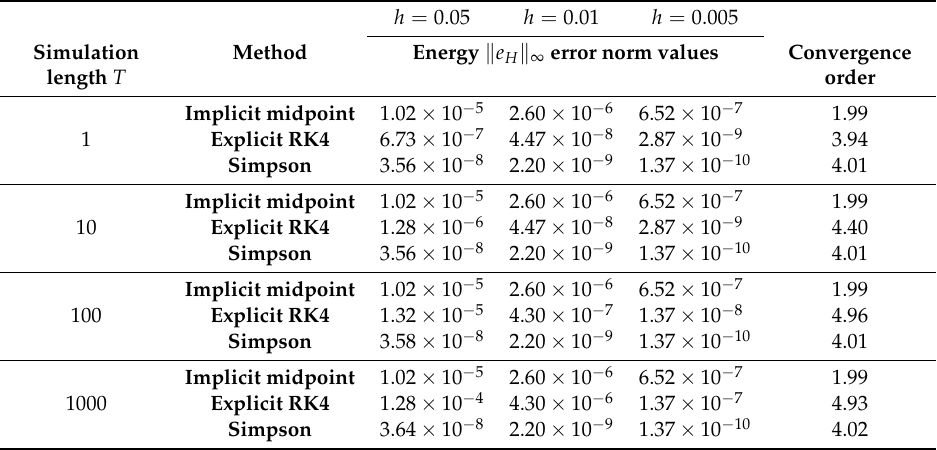}
	\end{adjustwidth}
\end{table}

\subsection{Comparison with a reduced Lagrange top system}
\label{subsec:reducedModel}

It is possible to reduce the second-order system \eqref{eq:ELtop}--\eqref{eq:conservedMomenta} to a first-order system by choosing special initial conditions, and by expressing it in terms of the conserved momenta $\{p_\varphi,p_\psi\}$. Let us choose the following  conditions:

\begin{itemize}
	\item inertias $I_3=2I=2$;
	\item initial nutation is set to $\theta_0=\sfrac{\pi}{6}$;
	\item initial precession and rotation are null, $\{\varphi_0=0,\psi_0=0\}$; 
	\item initial precession and nutation velocities are null, $\{\dot\varphi_0=0,\dot\theta_0=0\}$;
	\item by the conserved momentum $p_\psi$ (equations \eqref{eq:pPsi}--\eqref{eq:ELtop3}), and by the selected conditions above, $\dot\psi+\dot\varphi\cos\theta = \dot\psi_0$ is constant and chosen such that the initial rotation velocity is 1 round per second, $\dot\psi_0=2\pi$.
\end{itemize}

By system \eqref{eq:conservedMomenta} and the above conditions, the conserved momenta are
\begin{equation}
	p_\varphi = 2\pi\sqrt{3}\,;\qquad p_\psi = 4\pi\,.
\label{eq:pPhiReduced}
\end{equation} 

Equation \eqref{eq:ELtop3} can be written as
\begin{equation}
	\ddot\theta = \left( \frac{m_\mathrm{top}\cg\ell}{I} + \dot\varphi^2\cos\theta - \frac{p_\psi}{I}\dot\varphi \right)\sin\theta\,.
\label{eq:nutationSecondOrder}
\end{equation}
This equation can be stated as a first-order system by introducing an auxiliary variable $\xi$ such that
\begin{equation}
	\dot\theta=\xi\,,\qquad \dot\xi = \left( A + \dot\varphi^2\cos\theta - \frac{p_\psi}{I}\dot\varphi \right)\sin\theta.
\label{eq:nutationFirstOrder}
\end{equation}
where $A$ is a squared frequency, chosen such that the time constant $\bar t=1$: 
\begin{equation}
	A=\frac{(2\pi)^2}{\bar t}=4\pi^2.
\label{eq:frequencyA}
\end{equation}
The nutation evolution system \eqref{eq:nutationFirstOrder} is then coupled to the precession evolution. By system \eqref{eq:conservedMomenta}, $p_\varphi - p_\psi \cos\theta = I\dot\varphi\sin^2\theta$ and the precession evolution is given by
\begin{equation}
	\dot\varphi = \frac{p_\varphi-p_\psi\cos\theta}{I\sin^2\theta}\,.
\label{eq:precessionFirstOrder}
\end{equation}
Table \ref{tab:HSTparametersReduced} lists the special conditions used in this section.

\begin{table}[H] 
\sisetup{
	scientific-notation = true,
	round-mode = places,
	round-precision = 2}
\caption{Parameters for a reduced Lagrange top model and a motion with cusps.\label{tab:HSTparametersReduced}}
\newcolumntype{C}{>{\centering\arraybackslash}X}
\begin{tabularx}{\textwidth}{p{3.3cm}CCC}
	\toprule
	& {\textbf{Parameter}}
	& \textbf{Value}
	& \textbf{Unit}
	\\
	\midrule
	\multirow{4}{=}{Physical parameters of the top}
	& $m_\mathrm{top}$
	& 1
	& [\si{\kilogram}]
	\\
%	Cross-sectional moment of inertia 
	& $I$
	& 1
	& [\si{\kilogram\meter\squared}]
	\\
%	Longitudinal moment of inertia
	& $I_3$
	& 2
	& [\si{\kilogram\meter\squared}]
	\\
%	Center of mass location
	& $A$
	& $4\pi^2$
	& [\si{\per\second\squared}]
	\\
%	Center of mass location
	& $\ell$
	& $\sfrac{4\pi^2}{\cg}$ \textsuperscript{$\star$}
	& [\si{\meter}]
	\\
	\midrule
	\multirow{3}{*}{Initial conditions}
	& $\left\{\varphi_0,\theta_0,\psi_0\right\}$
	& $\left\{0,\,\sfrac{\pi}{6},\,0\right\}$
	& [\si{\radian}]
	\\[0.5ex]
	& $\left\{\dot\varphi_0,\dot\theta_0,\dot\psi_0\right\}$
	& $\left\{0,\,0,\,2\pi\right\}$
	& [\si{\radian\per\second}]
	\\
	\midrule
	\multirow{2}{*}{Conserved momenta}
	& $p_\varphi$
	& $2\pi\sqrt{3}$
	& [\si{\kilo\gram\meter\squared\per\second}]
	\\[0.5ex]
	& $p_\psi$
	& $4\pi$
	& [\si{\kilo\gram\meter\squared\per\second}]
	\\
	\bottomrule
\end{tabularx}
\noindent{\footnotesize{\textsuperscript{$\star$} $\cg$ is the gravity constant. Its value does not affect the results of this section.}\hfill}
\end{table}

The first order, \emph{reduced Lagrange top system}, is formed by equations \eqref{eq:nutationFirstOrder}--\eqref{eq:precessionFirstOrder}. Using the special conditions of Table \ref{tab:HSTparametersReduced}, the reduced system becomes
\begin{subequations}
	\begin{align}
		\label{eq:reducedSystemPhi}
		\dot\varphi
		& = 2\pi \left(\frac{\sqrt{3} - 2\cos\theta}{\sin^2\theta}\right)
		\\
		\label{eq:reducedSystemTheta}
		\dot\theta
		& = \xi
		\\
		\label{eq:reducedSystemXi}
		\dot\xi
		& = \left(4\pi^2 + \dot\varphi^2 \cos\theta - 4\pi \dot\varphi\right)\sin\theta\,.
	\end{align}
	\label{eq:reducedSystem}
\end{subequations}
\subsubsection{Motion with cusps and the reduced Lagrange top system}

\begin{figure}[H]
	\begin{adjustwidth}{-\extralength}{-\extralength}
	\centering
		\includegraphics[]{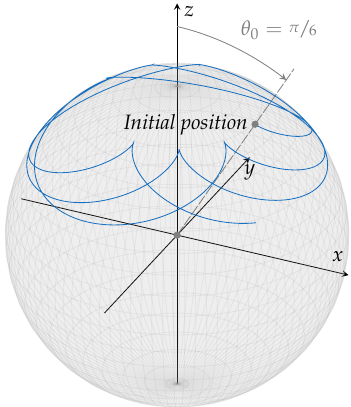}
		\qquad
		\includegraphics[]{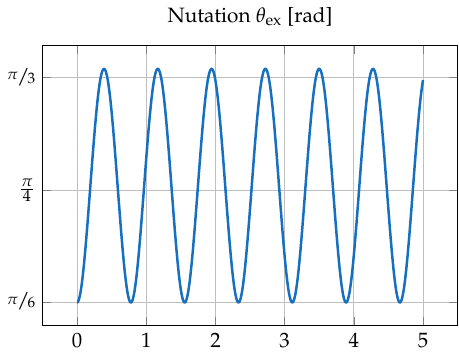}
	\end{adjustwidth}
	\caption{Motion of the reduced Lagrange top model. This particular example showcases cusps during motion.
	}
	\label{fig:topMotion3Dbis}
\end{figure}

\begin{figure}[H]
	\begin{adjustwidth}{-\extralength}{-\extralength}
	\centering
		\includegraphics[]{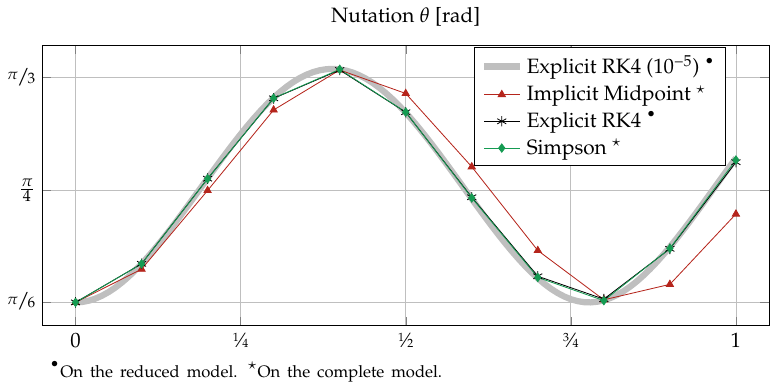}
	\end{adjustwidth}
	\caption{Several methods compared against the reference solution over one second for $h=0.1$. The reference solution is given by the RK4 method applied to the reduced model (equations \eqref{eq:reducedSystem}) of the heavy symmetric top motion at $h=10^{-5}$. The Simpson integrator is applied
to the complete model (equations \eqref{eq:ELtop}--\eqref{eq:conservedMomenta})	and correctly follows the reference solution.}
	\label{fig:nutationBis}
\end{figure}

In the rest of this section, the numerical approximations of the complete Lagrange top system \eqref{eq:ELtop}--\eqref{eq:conservedMomenta} are provided by the proposed Simpson integrator. Conversely, the numerical approximations of the reduced Lagrange top system \eqref{eq:reducedSystem} are provided by the RK4 method. The conditions of Table \ref{tab:HSTparametersReduced} are taken for these trials. Figure \ref{fig:topMotion3Dbis} shows the simulated motion showcasing cusps. Figure \ref{fig:nutationBis} shows the solutions obtained by the three integrators (Implicit midpoint, RK4, and Simpson) on this example. 
Even if RK4 solves the reduced system, Simpson tracks the reference solution with more precision for an equal step size.

\subsubsection{Further analysis on convergence}
\label{subsubsec:convergenceBis}

The reduced system \eqref{eq:reducedSystem} is of the first order but is different than the complete one  \eqref{eq:ELtop}--\eqref{eq:conservedMomenta}. The idea now is to provide another way to analyze convergence, in the case where an exact solution is not available. It is well-established that the RK4 method is fourth-order convergent. Therefore, if the two methods (Simpson and RK4) give the same result on different systems (complete Lagrange top and reduced Lagrange top), up to some error, their difference shall be consistent with the established convergence order of RK4.

\begin{figure}[H]
	\begin{adjustwidth}{-\extralength}{-\extralength}
	\centering
		\includegraphics[]{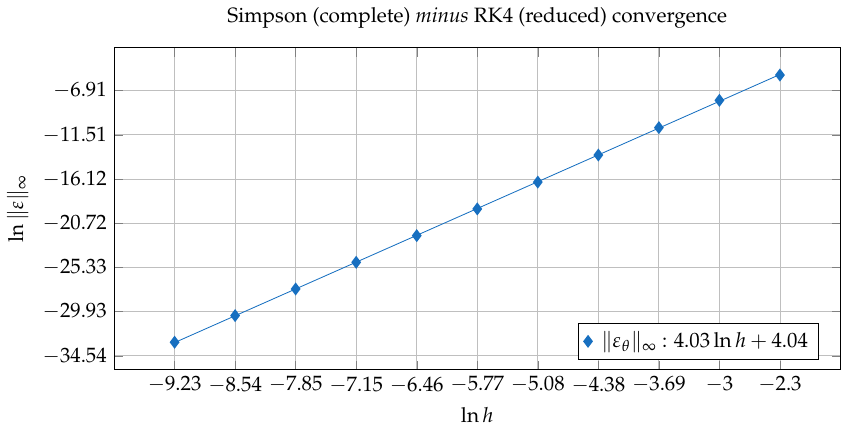}
	\end{adjustwidth}
	\caption{Consistency in the error convergence rate of the Simpson integrator on the full system, towards the RK4 solution of the reduced system. This difference is fourth-order convergent.}
	\label{fig:thetaConvergenceReduced}
\end{figure}

Let us calculate the $\ell^\infty$ error norm between the solution given by the Simpson integrator ($\theta_{\text{SIM}_n}$) on the complete system \eqref{eq:ELtop}--\eqref{eq:conservedMomenta}, and the solution provided by RK4 ($\theta_{\text{RK4}_n}$) on the reduced system \eqref{eq:reducedSystem}, by  evaluating 
\begin{equation*}
	\|\varepsilon_\theta\|_\infty 
	= 
	\sup_n \left\lvert \theta_{\text{SIM}_n}-\theta_{\text{RK4}_n} 
	\right\lvert
\end{equation*}
on diminishing step sizes $h$ successively divided by 2. The ratio between successive errors gives a number $2^{-\alpha}$ where $\alpha$ is a real number. The convergence order of this difference corresponds to the closest integer to $\alpha$. Figure \ref{fig:thetaConvergenceReduced} confirms the results that were previously obtained, when the Simpson approximation was compared with the exact solution of the Lagrange top. Unlike previous tests, this trial does not prove convergence but gives a very good indication of the convergence rate in the absence of an exact solution. The convergence order of the difference Simpson (complete) \emph{minus} RK4 (reduced) is consistent with the convergence order of the RK4 method, as expected.

The results of this section illustrate the precision and accuracy of Simpson's variational integrator. Let us remark that the convergence order was first analyzed by taking the analytical solution for the nutation of the Lagrange top, then taking a numerical solution on a reduced model with a fourth-order method. The results of this section indicate that the proposed Simpson integrator:
\begin{itemize}
	\item	does not artificially dissipate energy (Figure \ref{fig:energyTop});
	\item	exactly preserves the conserved momenta (Figure \ref{fig:MomentumConservation});
	\item	is fourth-order convergent on nutation and energy (Figures \ref{fig:thetaConvergenceREL} and \ref{fig:thetaConvergenceReduced});
	\item 	preserves its convergence rate for long simulations (Tables \ref{tab:nutationErrorNorms} and \ref{tab:energyErrorNorms}).
\end{itemize} 
\section{Concluding remarks}
\label{sec:conclusion}

The proposed Simpson variational integrator is a special case of the Galerkin-type variational integrators that use a Lobatto quadrature and Lagrange polynomials for interpolation. This paper has presented this integrator, specifically formulated for nonlinear systems characterized by inseparable Hamiltonians. It has been evaluated on two nonlinear examples: the double pendulum and the Lagrange top. Numerical trials show that the proposed Simpson integrator is fourth-order convergent on energy and on the exact solution of the nutation of a spinning top. It exactly preserves the conserved momenta. This symplectic integrator does not artificially dissipate energy and preserves its convergence rate for long simulations. This work shows that the proposed integrator is effective for nonlinear systems. It should be of interest in disciplines requiring accuracy and precision, but low to moderate complexity, such as the control of nonlinear systems.

A nonlinear formulation of higher-order Lobatto variational integrators could be of interest in areas where more precision is prioritized. This formulation will be derived for the sixth-order integrator exposed in \cite{ober2014,dubois2025_GSI} in the near future. Motivated by the results obtained in \cite{tsang2015,capobianco2024}, which develop specialized variational integrators for dissipative systems, the higher-order Lobatto variational integrator will be formulated for the optimal control of nonlinear systems subject to dissipation to resolve some of the issues exposed in \cite{rojas2024mathematics}. Additionally, as suggested in \cite{hairer2006}, Gaussian quadratures may be effective for systems that are subject to numerical stiffness. Variational integrators based on these quadratures shall be explored in future developments. Last but not least, it is well known that using Euler angles to model spinning tops is not ideal. These are affected by singularities, and group formulations such as quaternions are preferred. A geometric integrator for the special orthogonal and Euclidean groups shall be investigated in future developments.

%%%%%%%%%%%%%%%%%%%%%%%%%%%%%%%%%%%%%%%%%

\section*{Acknowledgments} This research was initiated while the corresponding author held a one-year visiting position at the Laboratoire des Sciences du Numérique de Nantes (LS2N) of the Institut Mines-Télécom Atlantique (IMTA) of Nantes, France, in 2024; under the authorization of Secihti. The corresponding author acknowledges the support received from the IMTA and the Secihti.

\appendix

\section{Simpson scheme for nonlinear Lagrangian systems}
\label{app:simpson}

Taking the difference
\begin{equation*}
\begin{aligned}
%	p_\gamma - p_{j,\gamma} 
	p_{\gamma,j+1} - p_{\gamma,j}
	%& = p_{r,\gamma,j+1} + p_{\ell,\gamma,j+1}
	%\\
	& = \frac{h}{12} \partial_\gamma M_{\alpha\beta}( q_j )\,g_j^\alpha g_j^\beta 
		+ \frac{h}{3} \partial_\gamma M_{\alpha\beta}( q_{j+\sfrac12} )\,g_{j+\sfrac12}^\alpha g_{j+\sfrac12}^\beta
		+ \frac{h}{12} \partial_\gamma M_{\alpha\beta}(q_{j+1})\,g_{j+1}^\alpha g_{j+1}^\beta 
		\\
		& \quad\ 
		- \frac h6 \partial_\gamma V( q_j )
		- \frac{2h}3 \partial_\gamma V( q_{j+\sfrac12} )
		- \frac h6 \partial_\gamma V(q_{j+1})\, ,
\end{aligned}
\end{equation*}
and the sum
\begin{equation*}
\begin{aligned}
	p_{\gamma,j+1} + p_{\gamma,j}
	& = \frac23 M_{\gamma\beta}( q_j )\,g_j^\beta
		 + \frac43 M_{\gamma\beta}( q_{j+\sfrac12} )\,g_{j+\sfrac12}^\beta
		 + \frac23 M_{\gamma\beta}(q_{j+1})\,g_{j+1}^\beta
		 \\
		 & \quad\ 
		 - \frac{h}{12} \partial_\gamma M_{\alpha\beta}( q_j )\,g_j^\alpha g_j^\beta
		 + \frac{h}{12} \partial_\gamma M_{\alpha\beta}(q_{j+1})\,g_{j+1}^\alpha g_{j+1}^\beta
		 + \frac h6 \partial_\gamma V( q_j )
		 - \frac h6 \partial_\gamma V(q_{j+1}),
\end{aligned}
\end{equation*}
results in a system of equations describing the evolution of the system in the discrete space: 
\begin{subequations}
\label{eq:discSimpEL}
\begin{align}
% PRIMERA ECUACIÓN
	M_{\gamma\beta}(q)\,g^\beta
		- M_{\gamma\beta}( q_j )\,g_j^\beta
		- \frac{h}{2} 
		\partial_\gamma M_{\alpha\beta}( q_{j+\sfrac12} )
		\,g_{j+\sfrac12}^\alpha g_{j+\sfrac12}^\beta
		+ h \partial_\gamma V( q_{j+\sfrac12} )
	& = 0,
	\label{eq:discSimpEL1}
	\\[1ex]
%	
% SEGUNDA ECUACIÓN CON SPLIT
	\begin{split}
	 p_{\gamma} - p_{\gamma,j} 
		+ \frac h6 \partial_\gamma V( q_j )
		+ \frac{2h}{3}\partial_\gamma V( q_{j+\sfrac12} )
		+ \frac h6 \partial_\gamma V(q)
%	agregamos un espacio horizontal para mover la primera línea de la ecuación hacia la izquierda	
	\hphantom{abcdefghijklmnopqrs}
	\\
	- \frac{h}{12} \partial_\gamma M_{\alpha\beta}( q_j )\,g_j^\alpha g_j^\beta 
		- \frac{h}{3} \partial_\gamma M_{\alpha\beta}( q_{j+\sfrac12} )\,g_{j+\sfrac12}^\alpha g_{j+\sfrac12}^\beta
		- \frac{h}{12} \partial_\gamma M_{\alpha\beta}(q)\,g^\alpha g^\beta 
		& = 0,
	\end{split}	
	\label{eq:discSimpEL2}
	\\[1ex]
%
% TERCERA ECUACIÓN CON SPLIT
	\begin{split}	
	\frac16 M_{\gamma\beta}( q_j )\,g_j^\beta
		 + \frac23 M_{\gamma\beta}( q_{j+\sfrac12} )\,g_{j+\sfrac12}^\beta
		 + \frac16 M_{\gamma\beta}(q)\,g^\beta
		-\frac{p_{\gamma} + p_{\gamma,j}}{2}
%	agregamos un espacio horizontal para mover la primera línea de la ecuación hacia la izquierda	
	\hphantom{abcdefgh}
	\\
	- \frac{h}{12} \partial_\gamma M_{\alpha\beta}( q_j )\,g_j^\alpha g_j^\beta
		 + \frac{h}{12} \partial_\gamma M_{\alpha\beta}(q)\,g^\alpha g_{j+1}^\beta
		 + \frac{h}6 \partial_\gamma V( q_j )
		 - \frac{h}6 \partial_\gamma V(q)
	& = 0.
	\end{split}
	\label{eq:discSimpEL3}
\end{align}
\end{subequations}
Subscripts $j+1$ are dropped in the above equations by establishing that $q\equiv q_{j+1}$ and $p\equiv p_{j+1}$. Let us remark that Simpson's quadrature can be recognized in equations \eqref{eq:discSimpEL2} and \eqref{eq:discSimpEL3}. As such, equation \eqref{eq:discSimpEL2} is the discrete analogue of $\dot p = - \nabla V + \frac 12 \dot q\tran \nabla M(q) \dot q$. In the same manner, equation \eqref{eq:discSimpEL3} is the discrete analogue  $M(q)\dot q = p$ when $h\rightarrow 0$. These equations verify \eqref{eq:multibodyDynamics} as expected.

Note that the scheme \eqref{eq:discSimpEL} is implicit, and it is not possible to eliminate the middle point as in \cite{rojas2024axioms}. Equation \eqref{eq:discSimpEL1} corresponds to the evolution equation for the point in the middle of the interval $q_m$. Equations \eqref{eq:discSimpEL2} and \eqref{eq:discSimpEL3} are the equations for $p$ and $q$, respectively. Let us establish that
\begin{equation*}
\begin{aligned}
	F_{q_m,j}(q_m,p,q)
	& \equiv h M_{\gamma\beta}(q)\,g^\beta
		-h M_{\gamma\beta}( q_\ell )\,g_\ell^\beta
		+ h^2 \partial_\gamma V( q_m )\,
		- \frac{h^2}{2} \partial_\gamma M_{\alpha\beta}( q_m )\,g_m^\alpha g_m^\beta\,,
	\\
	F_{p,j}(q_m,p,q)
	& \equiv
	h \left(p_\gamma - p_{j,\gamma}\right) 
		+ h^2\left[ \frac16 \partial_\gamma V( q_\ell )
		+ \frac23 \partial_\gamma V( q_m )
		+ \frac16 \partial_\gamma V(q)
		\right]
	\\
	& \qquad\ 		
	- \frac{h^2}{2}\left[ \frac{1}{6} \partial_\gamma M_{\alpha\beta}( q_\ell )\,g_\ell^\alpha g_\ell^\beta 
		+ \frac{2}{3} \partial_\gamma M_{\alpha\beta}( q_m )\,g_m^\alpha g_m^\beta
		+ \frac{1}{6} \partial_\gamma M_{\alpha\beta}(q)\,g^\alpha g^\beta 
		\right],
	\\
	F_{q,j}(q_m,p,q)
	& \equiv
	h\left[ \frac16 M_{\gamma\beta}( q_\ell )\,g_\ell^\beta
		 + \frac23 M_{\gamma\beta}( q_m )\,g_m^\beta
		 + \frac16 M_{\gamma\beta}(q)\,g^\beta
	\right]
	-\frac{h}{2} \left(p_\gamma + p_{j,\gamma}\right)
	\\
	& \qquad\
	- \frac{h^2}{12} \partial_\gamma M_{\alpha\beta}( q_\ell )\,g_\ell^\alpha g_\ell^\beta
		 + \frac{h^2}{12} \partial_\gamma M_{\alpha\beta}(q)\,g^\alpha g^\beta
		 + \frac{h^2}6 \partial_\gamma V( q_\ell )
		 - \frac{h^2}6 \partial_\gamma V(q)\,.
\end{aligned}
\end{equation*}
For a fixed triple $(q_{j+\sfrac12},p_j,q_j)$, it is required that $F(q_m,p,q)=0$. Let us consider Newton's algorithm $(q_m,p,q)\rightarrow(\tilde q_m,\tilde p,\tilde q)$ where $\tilde q_m = q_m + \delta q_m$, $\tilde p = p +\delta p$ and $\tilde q = q + \delta q$:
\begin{equation*}
	F(q_m,p,q) + \diffp{F}{q_m} \delta q_m + \diffp Fp \delta p + \diffp Fq \delta q= 0,
\end{equation*}
by initializing $(q_m,p,q)=(q_j,p_j,q_j)$ at each iteration. Using $K_{\gamma\delta}(q_m) \equiv \partial_\gamma \partial_\delta V(q_m)$, the elements of the Jacobian are block matrices given by
\begin{equation*}
\begin{aligned}
	\diffp{F_{q_m,\gamma}}{q_m^\delta} 
	& = 
	-4\left[M_{\gamma\delta}(q) + M_{\gamma\delta}(q_\ell)\right] 
	+ h^2 K_{\gamma\delta}(q_m) 
	- \frac{h^2}{2} \partial_\gamma \partial_\delta M_{\alpha\beta}(q_m)g_m^\alpha g_m^\beta\,
	\\[0.5ex]
	\diffp{F_{q_m,\gamma}}{p_\delta} 
	& = 0\,
	\\[0.5ex]
	\diffp{F_{q_m,\gamma}}{q^\delta} 
	& = M_{\gamma\delta}(q_\ell) 
		+ 3M_{\gamma\delta}(q)
		- h \partial_\gamma M_{\alpha\delta}(q_m)g_m^\alpha 
		+ h \partial_\delta M_{\gamma\beta}(q_m)g_m^\beta\,
	\\[1.5ex]
	\diffp{F_{p,\gamma}}{q_m^\delta} 
	& = 
	-\frac{2h}3 \partial_\gamma M_{\alpha\delta}(q_\ell)g_\ell^\alpha
	+\frac{2h}3 \partial_\gamma M_{\alpha\delta}(q)g^\alpha
	-\frac{h^2}3 \partial_\delta \partial_\gamma M_{\alpha\beta}(q_m)g_m^\alpha g_m^\beta
	+\frac{2h^2}3 K_{\delta\gamma}(q_m)\,
	\\[0.5ex]
	\diffp{F_{p,\gamma}}{p_\delta} 
	& = h\delta_{\gamma\delta}\,
	\\[0.5ex]
	\diffp{F_{p,\gamma}}{q^\delta} 
	& = \frac{h}6 \partial_\gamma M_{\alpha\delta}(q_\ell)g_\ell^\alpha
	-\frac{2h}3 \partial_\gamma M_{\alpha\delta}(q_m)g_m^\alpha
	-\frac{h}2 \partial_\gamma M_{\alpha\delta}(q)g^\alpha
	\\
	& \qquad\
	-\frac{h^2}{12} \partial_\delta \partial_\gamma M_{\alpha\beta}(q)g^\alpha g^\beta
	+\frac{h^2}6 K_{\delta\gamma}(q)\,
	\\[1.5ex]
	\diffp{F_{q,\gamma}}{q_m^\delta} 
	& = \frac{2}{3} \left[
		M_{\gamma\delta}(q_\ell)
		- M_{\gamma\delta}(q)
	  	\right]
		+\frac{2h}3 \partial_\delta M_{\gamma\beta}(q_m)g_m^\beta
		-\frac{h}3 \partial_\delta M_{\alpha\delta}(q_\ell)g_\ell^\beta
		-\frac{h}3 \partial_\delta M_{\alpha\delta}(q)g^\beta\,
	\\[1ex]
	\diffp{F_{q,\gamma}}{p_\delta} 
	& = -\frac12 \delta_{\gamma\delta}\,
	\\[1ex]
	\diffp{F_{q,\gamma}}{q^\delta} 
	& = - \frac16 M_{\gamma\delta}(q_\ell)
		+ \frac23 M_{\gamma\delta}(q_m)
		+ \frac12 M_{\gamma\delta}(q)
		+ \frac{h^2}{24}\partial_\delta \partial_\gamma M_{\alpha\beta}(q)g^\alpha g^\beta
		- \frac{h^2}{12} K_{\gamma\delta}
		\\
		& \qquad\
		+ \frac{h}6 \partial_\delta M_{\gamma\beta}(q)g^\beta
		+ \frac{h}{12} \partial_\gamma M_{\alpha\delta}(q_\ell)g_\ell^\beta
		+ \frac{h}4 \partial_\delta M_{\alpha\delta}(q)g^\beta\,.
\end{aligned}
\end{equation*}

\section{Implicit midpoint integrator for nonlinear systems}
\label{app:Newmark}

The Implicit midpoint scheme \cite{simo1992} is a classical variational integrator that is very popular due to its simplicity. It is usually taken as a benchmark variational integrator, and its properties have been given in \cite{simo1992,sanzSerna1992,marsden2001}. For linear problems, the Implicit midpoint method is equivalent to Newmark's method \cite{newmark1959}, which is very popular for structural dynamics problems \cite{geradin2015,chopra2020}. In this section, we provide its nonlinear formulation using similar notations to those of the proposed Simpson's integrator.

Let us take the subscripts $(\ell,m)$ implying discrete ``left'' and ``middle'' values, respectively; ``right'' values do not use a subscript. Let us consider the corresponding positions $(q_\ell,q_m,q)$. The states $q_\ell$ and $q$ are separated by a time step $h\in\mathbb{R}$ so that for an interval $[0,h]$, $q_\ell\approx q(0)$ and $q_r\approx q(h)$. The state at the middle of the discretization interval is
\begin{equation*}
	q_m^\alpha = \frac12 \left(q_\ell^\alpha + q^\alpha \right)
\end{equation*}
The implicit midpoint discretization considers a finite difference $g\in Q\times Q$ as the discrete analogue of a velocity vector: 
\begin{equation*}
	g^\alpha = \frac1h\left(q^\alpha - q_\ell^\alpha \right).
\end{equation*}
The discrete Lagrangian then approximates the action integral $\mathcal S$ along the curve segment between $q_\ell$ and $q$. The approximation to $\int_0^T \mathcal L\dl t$ is given by the middle point quadrature:
\begin{equation*}
	L_\mathrm{d} = \frac h2 g^\gamma M_{\gamma\beta}(q_m)g^\beta - h V (q_m).
\end{equation*}
Considering a discrete curve of points $\{q_0,\ldots,q_N\}$, the action sum is
\begin{equation*}
	S_\mathrm{d} = \ldots + L_\mathrm{d}\left(q_{j-1},q_j\right) + L_\mathrm{d}\left(q_j,q_{j+1}\right) + \ldots.
\end{equation*}
The variational principle $\delta  S_d =0$ for an arbitrary $\delta q_j$ yields the discrete Euler-Lagrange equations
\begin{equation*}
	\diffp{L_d}{q_r}\left(q_{j-1},q_j\right) + \diffp{L_d}{q_\ell}\left(q_j,q_{j+1}\right) = 0.
\end{equation*}
The generalized momentum is defined on the right as $p_{j}=\diffp{L_d}{q_r}\left(q_{j-1},q_j\right)$. As such, the left term of the above Euler-Lagrange equations can be identified as being the discrete generalized momentum, and so 
\begin{equation*}
	p_{\alpha,j} 
	= -\diffp{L_d}{q_\ell}\left(q_j,q_{j+1}\right)\,;
	\qquad
	p_{\alpha,j+1} 
	= -\diffp{L_d}{q_r}\left(q_j,q_{j+1}\right)\,.
\end{equation*}
Taking the difference between these two expressions yields
\begin{equation*}
	p_{\alpha,j+1} - p_{\alpha,j}
	= \frac h2 \left(\partial_\alpha M_{\gamma\beta}\right)g^\gamma g^\beta - h\left(\partial_\alpha V\right)\,,
\end{equation*}
and taking the sum yields
\begin{equation*}
	p_{\alpha,j+1} + p_{\alpha,j}
	= 2M_{\alpha\beta} g^\beta\,
\end{equation*}
results in a system of equations describing the evolution of the system in the discrete phase space:
\begin{subequations}
	\label{eq:NewmDiscreteMotionEq}
	\begin{align}
		\frac{p_\alpha - p_{\alpha,j} }{h}
		- \frac 12 \left(\partial_\alpha M_{\gamma\beta}\right)g^\gamma g^\beta 
		+ \left(\partial_\alpha V\right)
		& = 0
		\\
		M_{\alpha\beta}\, \frac{q^\beta - q_j^\beta}{h} 
		- \frac 12 \left(p_\alpha + p_{\alpha,j}\right)
		& = 0\,.
	\end{align}
\end{subequations}
In the above, we have suppressed the subscript ``$j+\sfrac 12$'' since it is not required anymore. Note that equations \eqref{eq:NewmDiscreteMotionEq} are discrete analogues of $\dot p = - \nabla V + \frac 12 \dot q\tran \nabla M(q) \dot q$ and $M(q)\dot q = p$, which verify equations \eqref{eq:multibodyDynamics} as expected.

The scheme \eqref{eq:NewmDiscreteMotionEq} is implicit. Let us establish that
\begin{equation*}
	F(p,q) = 
	\def\arraystretch{1.5}
	\begin{pmatrix}
		F_{p,j}(p,q)
		\\
		F_{q,j}(p,q)
	\end{pmatrix}
	\equiv 
	\begin{pmatrix}
		p_\alpha - p_{\alpha,j} - \frac h2 \left(\partial_\alpha M_{\gamma\beta}\right)g^\gamma g^\beta + h\left(\partial_\alpha V\right)
		\\
		- \frac h2 \left(p_\alpha + p_{\alpha,j}\right)
	\end{pmatrix}\,.
\end{equation*}
For a fixed tuple $\left(p_j,q_j\right)$, it is required that $F(p,q) = 0$. Let us consider Newton's algorithm $(p,q)\rightarrow(\tilde p,\tilde q)$ where $\tilde p = p +\delta p$ and $\tilde q = q + \delta q$:
\begin{equation*}
	F(p,q) + \diffp Fp \delta p + \diffp Fq \delta q = 0,
\end{equation*}
by initializing $(p,q)=(p_j,q_j)$ at each iteration. Then, the elements of the Jacobian are given by
\begin{equation*}
	\begin{pmatrix}
		\diffp{F_{p,j,\alpha}}{p_\mu}
		&
		\diffp{F_{p,j,\alpha}}{q^\mu}
		\\
		\diffp{F_{q,j,\alpha}}{p_\mu}
		&
		\diffp{F_{q,j,\alpha}}{q^\mu}
	\end{pmatrix}
	= 
	\begin{pmatrix}
		\delta_\alpha^\mu
		&
		\frac h2 K_{\alpha\mu} - \frac h4\left(\partial_{\alpha\mu}^2 M_{\gamma\beta}\right) g^\gamma g^\beta - \left(\partial_\alpha M_{\mu\beta}\right) g^\beta
		\\
		-\frac h2 \delta_\alpha^\mu
		&
		M_{\alpha\mu} + \frac h2 \left(\partial_\mu M_{\alpha\beta}\right) g^\beta
	\end{pmatrix}
\end{equation*}
with $K_{\alpha\mu}(q_m) \equiv \left(\partial_{\alpha\mu}^2 V\right)\!\!(q_m)$.

\bibliographystyle{unsrturl}
\bibliography{NLSimpsonTop.bib}

\end{document}